\documentclass{article}

\usepackage{amsfonts}
\usepackage{amsmath}
\usepackage{amssymb}
\usepackage{amsthm}
\usepackage{amscd}
\usepackage{verbatim}
\usepackage[dvips]{graphicx}
\usepackage{graphics}
\usepackage[all,2cell]{xypic}

\theoremstyle{plain}
\newtheorem{thm}{Theorem}[section]

\newtheorem{cor}[thm]{Corollary}
\newtheorem{lem}[thm]{Lemma}
\newtheorem{prop}[thm]{Proposition}

\theoremstyle{definition}

\newtheorem{ex}[thm]{Example}
\newtheorem{nt}[thm]{Notations}
\newtheorem{rem}[thm]{Remark}

\theoremstyle{remark}

\DeclareMathOperator{\Spec}{Spec}

\DeclareMathOperator{\Cone}{Cone}
\DeclareMathOperator{\Cyl}{Cyl}
\DeclareMathOperator{\im}{Im}
\DeclareMathOperator{\Ker}{Ker}

\DeclareMathOperator{\Ar}{Ar}
\DeclareMathOperator{\id}{id}

\DeclareMathOperator{\E}{E}
\DeclareMathOperator{\Homo}{H}

\DeclareMathOperator{\AAA}{\mathcal{A}}
\DeclareMathOperator{\BBB}{\mathcal{B}}
\DeclareMathOperator{\CCC}{\mathcal{C}}

\DeclareMathOperator{\FFF}{\mathcal{F}}

\DeclareMathOperator{\MMM}{\mathcal{M}}
\DeclareMathOperator{\PPP}{\mathcal{P}}
\DeclareMathOperator{\WWW}{\mathcal{W}}
\DeclareMathOperator{\XXX}{\mathcal{X}}
\DeclareMathOperator{\YYY}{\mathcal{Y}}
\DeclareMathOperator{\ZZZ}{\mathcal{Z}}

\DeclareMathOperator{\Kos}{\textbf{Kos}}

\DeclareMathOperator{\isoto}{\overset{\scriptstyle{\sim}}{\to}}

\DeclareMathOperator{\qis}{qis}
\DeclareMathOperator{\isom}{isom}

\UseAllTwocells

\title{Higher algebraic $K$-theory of finitely generated torsion modules over principal ideal domains}
\author{Satoshi Mochizuki
\footnote{This research is supported by the 21 century COE program at Graduate School of Mathematical Sciences, the University of Tokyo.
}}
\date{}

\begin{document}

\maketitle

\vspace{-1cm}
\begin{center}
\tt{mochi@ms.u-tokyo.ac.jp}
\end{center}

\begin{abstract}
The main purpose of this paper is computing higher algebraic $K$-theory of Koszul complexes over principal ideal domains. The second purpose of this paper is giving examples of comparison techniques on algebraic $K$-theory for Waldhausen categories without the factorization axiom.  
\end{abstract}


\tableofcontents

\section{Introduction}

\noindent
Let $R$ be a non field, principal ideal domain, $\MMM^1(R)$ the category of finitely generated torsion $R$-modules and $\Kos^1(R)$ the full subcategory of complexes of finitely generated $R$-modules such that $X_i=0$ for $i \ne 0,1$ and its boundary morphism $d^X:X_1 \to X_0$ is injective and whose homology groups are in $\MMM^1(R)$. Then $\Kos^1(R)$ and $\MMM^1(R)$ can be considered as exact categories in the natural way. The one of the purpose of this note, we will prove the following theorems.

\begin{thm}[Comparison theorem] \label{comparison} $ $ 
\\
An exact functor $\Homo_0:\Kos^1(R) \to \MMM^1(R)$ induces an isomorphism
$$K_n(\Kos^1(R);\qis) \isoto K_n(\MMM^1(R))$$
for any non-negative integer $n$.
\end{thm}

\begin{thm}[Split fibration theorem] \label{split fib} $ $
\\
In the notation above, we have a split short exact sequence
$$0 \to K_n(\Kos^1(R)^{q};\isom) \to K_n(\Kos^1(R);\isom) \to K_n(\Kos^1(R);\qis) \to 0$$
where $\Kos^1(R)^{q}$ is the full subcategory of acyclic complexes in $\Kos^1(R)$.
\end{thm}

\noindent
\textbf{Theorem \ref{comparison}} is related with the author's work about Gersten's conjecture. (See \cite{Moc07}.) The difficulty of verifying to the theorems above is lack of a cylinder functor on $\Kos^1(R)$. In general, several important theorems for Waldhausen $K$-theory need the hypothesis that an underlying Waldhausen category satisfies the factorization axioms. (See \cite{Sch06}). So an another purpose of this note, we will establish several techniques on algebraic $K$-theories for Waldhausen categories without the factorization axiom. To do so, in \S 2, we will examine the cell filtrations theorem in \cite{Wal85}. The hard part of driving this theorem is checking the hypothesis (=existence of cellular factorizations). So we will give an example of how to construct them in a certain case. In \S 3, we will ad hoc define a concept of semi-direct product of exact categories. Using this notion, we will propose an example of a technique on constructing split exact sequences only using the result \cite{Qui73} and \textbf{Theorem \ref{comparison}}. In appendix, we will give an another proof of \textbf{Theorem \ref{split fib}} by using Schlichting brilliant theorem in \cite{Sch04}.\\
\\
\noindent
\textbf{Prerequisite knowledge} In this paper, we will frequently used the theory in \cite{Qui73}, \cite{Wal85}, \cite{TT90}. (When we quote the paper \cite{Qui73}, we will follow the page numbering by the numbering at the bottom.)\\
\\
\noindent
\textbf{Acknowledgement} This works except appendix was done four years ago. The author thankful to Shuji Saito for recommending for him to write up the paper, to Masana Harada for training his techniques on algebraic $K$-theory, to Fabrice Orgogozo for kindly responding to his delusional strategies of how to extend Theorem \ref{comparison} to general regular domain, to Takeshi Saito for making him to get to the core of the proof Theorem \ref{comparison}, and to Toshiro Hiranouchi for carefully reading a draft version of this paper.

\section{Notations and remarks}

In this section, we will keep the notations, and give some remarks. First
notice that there are many choices of defining a cylinder functor, sign notations and truncation functors about chain complexes over an additive (or abelian) category. So we will keep the notations. 

\begin{nt} [General notations of chain complexes] \label{gncc} $ $
\\
(1) Let $\XXX$ be an additive category. Consider chain complexes $X$ in $\XXX$. 
We use the topologist's indexing, so differentials decreasing degree: 
$d_{n+1}^X:X_{n+1} \to X_n$.\\
We denote $C(\XXX)$ as the category of chain complexes in $\XXX$, and $C_b(\XXX)$ full sub-category of bounded complexes in $C(\XXX)$.
\\
(2) For any integer $k$ and complex $X$ in $\XXX$, we define complex $X[k]$ as $X[k]_n=X_{n+k}$ and $d_n^{X[k]}={(-1)}^kd_{n+k}^X$.\\
(3) For a morphism of complexes $f:X \to Y$ in $\XXX$, we define the complexes $\Cyl (f)$ and $\Cone (f)$ as the following way:
$$ \Cyl (f)_n =X_n \oplus X_{n-1} \oplus Y_n , \ \ d_n^{\Cyl f}=
\begin{pmatrix}
d_n^X & \id & 0\\
0 & -d_{n-1}^X & 0\\
0 & -f_{n-1} & -d_n^Y
\end{pmatrix}.$$
$$ \Cone (f)_n =X_{n-1} \oplus Y_n , \ \ d_n^{\Cone f}=
\begin{pmatrix}
-d_{n-1}^X & 0\\
-f_{n-1} & -d_n^Y
\end{pmatrix}.$$
(4) The operation $\Cyl$ is considered as the functor
$$\Cyl :\Ar C(\XXX) \longrightarrow C(\XXX)$$
where $\Ar C(\XXX)$ means the category of arrows in $C(\XXX)$.
That is, to each commutative square in $C(\XXX)$
$$\xymatrix{
X \ar[r]^{f} \ar[d]^{a} & Y \ar[d]_{b}\\
X' \ar[r]^{g} & Y'
}$$
$\Cyl$ assigns functorially a morphism $\Cyl(a,b):
\Cyl(f) \to \Cyl(g)$. This is defined by
$$ {\Cyl(a,b)}_n=
\begin{pmatrix}
a_n & 0 & 0 \\
0 & a_{n-1}& 0\\
0 & 0 & b_n
\end{pmatrix}.$$
(5) There are specific natural transformations $j_1$, $j_2$, and $p$
such that for each morphism $f:X \to Y$ in $C(\XXX)$, 
the following diagram is commutative:
$$\xymatrix{
X \ar[r]^{\scriptstyle{j_1(f)}} \ar[dr]_{f} & \Cyl(f) \ar[d]_{\scriptstyle{p(f)}} & Y \ar[l]_{\scriptstyle{j_2(f)}} \ar[ld]^{\id_{Y}}\\
& Y & .  
}$$
They are defined by
$$j_1(f)_n:=
\begin{pmatrix}
\id_{X_n} \\
0 \\
0
\end{pmatrix}
,\ \ 
j_2(f)_n:=
\begin{pmatrix}
0 \\
0 \\
\id_{Y_n}
\end{pmatrix}
,\ \ 
p(f)_n:=
\begin{pmatrix}
f_n & 0 & \id_{Y_n}
\end{pmatrix}.$$
It is well known that $p(f):\Cyl(f) \to Y$
is a homotopy equivalence. More precisely
$p(f) \circ j_2(f) = \id_{Y}$ and
$j_2(f) \circ p(f)$ is chain homotopic to 
$\id_{\Cyl(f)}$.\\ 
(6) If $\XXX$ be an abelian category, then for any chain complex $X$ in $\XXX$ and an integer $n$, we can define its {\it{truncation complexes}} $\tau_{\geqq n}X$ and $\tau_{\leqq n}X$ as follows:
$$ \tau_{\geqq n}X: \cdots \to X_{n+2} \to X_{n+1} \to \Ker d_n^X \to 0 \to 0 \to \cdots$$
$$ \tau_{\leqq n}X: \cdots \to 0 \to 0 \to \im d_{n+1}^X \to X_n \to X_{n-1} \to \cdots .$$  
So for any integer $n$, there is a canonical short exact sequence:
\begin{equation} \label{short exact}
0 \to \tau_{\geqq n+1}X \to X \to \tau_{\leqq n} X \to 0  \ .  
\end{equation}
\end{nt}

\noindent
Next we will change some notations in \cite{Wal85} into as follows.

\begin{nt}[$S_{\bullet}$-construction] $ $ 
\\
In this paper, for any Waldhausen category $(\XXX,w)$, we will change the notation,  $wS_{\bullet}\XXX$ in \cite{Wal85} into $S_{\bullet}(\XXX;w)$ and denote its $n$-th $K$-group by $K_n(\XXX;w)$. If $w=\isom$, that is, the class of all isomorphisms, we will abbreviate as $K_n(\XXX)$.  
\end{nt}

\noindent
The following statement may be well-known. But to accomplish logical completeness, we will give the proof. 

\begin{prop} \label{exact check} $ $
\\
Let $\XXX$ and $\YYY$ be small exact categories.\\
{\rm (1)} An exact functor $f:\XXX \to \YYY$ in the sense of Quillen {\rm \cite{Qui73}} p.100 is equivalent to an exact functor in the sense of Waldhausen {\rm \cite{Wal85}} p.327, {\rm \cite{TT90}} p.254 1.2.7.\\
{\rm (2)} An admissible exact sequence $f' \rightarrowtail f \twoheadrightarrow f''$ of exact functors $\XXX \to \YYY$ in the sense of Quillen {\rm \cite{Qui73}} p.106 is equivalent to a cofibration sequence of exact functors in the sense of Waldhausen {\rm \cite{Wal85}} p.331, {\rm \cite{TT90}} p.262 1.7.3.\\
{\rm (3)} The exact category of exact sequences in the sense of Quillen {\rm \cite{Qui73}} p.105, {\rm \cite{Wal78}} p.183 is coincided with those in the sense of Waldhausen {\rm \cite{Wal85}} p.328, {\rm \cite{TT90}} p.261 1.7.1.
\end{prop}

\begin{proof}
(1) The statement means that if $f$ preserves admissible short exact sequences, then $f$ preserves pushout along an admissible monomorphism, that is, the canonical map
$$fC \underset{fA}{\coprod} fB \to f(C \underset{A}{\coprod} B)$$
is an isomorphism whenever $A \rightarrowtail B$ is an admissible monomorphism. It is easily follows from \textbf{Lemma \ref{exact (co)base-change}} below.\\
(2) The statement means that if it has the property that for every $A$ in $\XXX$, the sequence
$$f'(A) \rightarrowtail f(A) \twoheadrightarrow f''(A)$$
is an admissible exact sequence in $\YYY$, then
it has also the property that for every admissible monomorphism 
$A \rightarrowtail B$ in $\XXX$ the square of admissible monomorphisms
$$\xymatrix{
F'(A) \ar@{>->}[r] \ar@{>->}[d] & F'(B) \ar@{>->}[d]\\
F(A) \ar@{>->}[r] & F(B)
}$$
is admissible in the sense that 
$F(A) \underset{F'(A)}{\coprod} F'(B) \rightarrowtail F(B)$
is also an admissible monomorphism. Let $A \rightarrowtail B \twoheadrightarrow C$ be an admissible exact sequence in $\XXX$ and consider the following commutative diagram
$$\xymatrix{
F'(A) \ar@{>->}[r] \ar@{>->}[d] & F(A) \ar@{->>}[r] \ar@{>->}[d] & F''(A) \ar@{>->}[d]\\
F'(B) \ar@{>->}[r] \ar@{->>}[d] & F(B) \ar@{->>}[r] \ar@{->>}[d] & F''(B) \ar@{->>}[d]\\
F'(C) \ar@{>->}[r] & F(C) \ar@{->>}[r] & F''(C) & .
}$$
Applying \textbf{Lemma \ref{Qui ex=Wal ex}} below for the diagram above, we get the result.\\
(3) Let $\ZZZ$ and $\WWW$ be sub exact categories and $\E(\ZZZ,\XXX,\WWW)$ be the exact category of exact sequences $Z \rightarrowtail X \twoheadrightarrow W$ such that $Z$ is in $\ZZZ$ and $W$ is in $\WWW$. The problem is that an admissible exact sequence in $\E(\ZZZ,\XXX,\WWW)$ is also a cofibration sequence. That is, we want to check that in the $3 \times 3$-diagram in the statement of \textbf{Lemma \ref{Qui ex=Wal ex}}, we have the admissible exact sequence (\ref{Qui ex=Wal ex 1}). So it is just a consequence of \textbf{Lemma \ref{Qui ex=Wal ex}}.
\end{proof}

\begin{lem} \label{exact (co)base-change} $ $
\\
Let $\XXX$ be a small exact category. Consider the diagram of admissible 
short exact sequences in $\XXX$
$$\xymatrix{
X  \ar@{>->}[r] \ar[d]_a \ar@{}[dr] |{1} & Y \ar@{->>}[r] \ar[d]^b \ar@{}[dr] |{2} & Z \ar[d]^c\\
X' \ar@{>->}[r] & Y' \ar@{->>}[r] & Z' & .  
}$$
Then the following are satisfied.\\
{\rm {(1)}} 1 is a coCartesian square if and only if $c$ is an isomorphism.\\
{\rm {(2)}} 2 is a Cartesian square if and only if $a$ is an isomorphism.
\end{lem}

\begin{proof}
By the Gabriel-Quillen embedding theorem (see \cite{Gab62}, \cite{Qui73} p.100, \cite{Kel90} p.408 A.2 
\cite{TT90} p.399 A.7.1.), we can reduce the problem to well-known statements in an abelian category.
\end{proof}

\begin{lem} \label{Qui ex=Wal ex} $ $
\\
Let $\XXX$ be a small exact category. Consider the diagram of admissible short exact sequences in $\XXX$
$$\xymatrix{
X \ar@{>->}[d]_f \ar@{>->}[r]^{i^X} & X' \ar@{>->}[d]_{f'} \ar@{->>}[r]^{p^X} & X'' \ar@{>->}[d]^{f''}\\
Y \ar@{->>}[d]_g \ar@{>->}[r]^{i^Y} & Y' \ar@{->>}[d]_{g'} \ar@{->>}[r]^{p^Y} & Y'' \ar@{->>}[d]^{g''}\\
Z \ar@{>->}[r]^{i^Z} & Z' \ar@{->>}[r]^{p^Z} & Z'' & .
}$$
Then the following two sequences are admissible exact sequences.
\begin{equation} \label{Qui ex=Wal ex 1}
Y \underset{X}{\coprod} X' \rightarrowtail Y' 
\overset{\scriptstyle{p^Z \circ g'}}{\twoheadrightarrow} Z''\ \ \ ,
\end{equation}
\begin{equation} \label{Qui ex=Wal ex 2}
X \overset{\scriptstyle{i^Y \circ f}}{\rightarrowtail} 
Y' \twoheadrightarrow Z'' \underset{Z''}{\times} Y'' \ \ \ ,
\end{equation}
where the morphism $Y \underset{X}{\coprod} X' \rightarrowtail Y'$ 
{\rm (}resp. $Y' \twoheadrightarrow Z'' \underset{Z''}{\times} Y''$ {\rm )} 
is induced from $i^Y$ and $f'$ 
{\rm (} resp. $p^Y$ and $g'$ {\rm )} 
by universal property of $Y \underset{X}{\coprod} X'$
{\rm (} resp. $Z'' \underset{Z''}{\times} Y''$ {\rm )}.
\end{lem}

\begin{proof}
This is essentially proved in \cite{TT90} p.262 1.7.4. We only prove the exactness of (\ref{Qui ex=Wal ex 1}). The exactness of (\ref{Qui ex=Wal ex 2}) is just the matter of dualizing the argument of that of (\ref{Qui ex=Wal ex 1}). First consider the following commutative diagram.
$$\xymatrix{
X \ar@{>->}[r]^{i^{X}} \ar@{>->}[d]_{f} & X' \ar@{->>}[r]^{p^X} \ar[d] & X'' \ar[d]^{\id_{X''}}\\
Y \ar@{>->}[r] & Y \underset{X}{\coprod} X' \ar[r] & X'' & ,
}$$  
where the morphism $Y \rightarrowtail Y \underset{X}{\coprod} X'$ is a co-base change of $i^X$ and
$Y \underset{X}{\coprod} X' \to  X''$ is induced from $p^X$ and zero morphism $Y \to X''$ by the universal property of $Y \underset{X}{\coprod} X'$.
Then by \textbf{Lemma \ref{exact (co)base-change}} (1), we learn the following sequence is admissible short exact sequence.
$$Y \rightarrowtail  Y \underset{X}{\coprod} X' \twoheadrightarrow  X'' \ \ \ .$$
Thus the diagram below has admissible exact sequences as rows
$$\xymatrix{
Y \ar@{>->}[r] \ar[d]_{\id_Y} & Y \underset{X}{\coprod} X' \ar[r] \ar[d] & X'' 
\ar@{>->}[d]_{f''}\\
Y \ar@{>->}[r]_{i^Y} & Y' \ar@{->>}[r]^{p^Y} \ar@{->>}[d]_{g'' \circ p^Y} 
& Y'' \ar@{->>}[d]_{g''}\\
& Z'' \ar[r]_{\id_{Z''}} & Z'' & . 
}$$
Then using the snake lemma (More precisely, we fix the Gabriel-Quillen embedding $\XXX \hookrightarrow \YYY$ and consider in $\YYY$), we learn the sequence 
(\ref{Qui ex=Wal ex 1}) is admissible exact.   
\end{proof}

\section{Examining the cell filtrations theorem}

In this section we will prove \textbf{Theorem \ref{comparison}}.

\begin{nt} $ $
\\
(1) Let $\AAA$ be the full additive subcategory of 
bounded complexes of finitely generated free $R$-modules with 
homology groups in $\MMM^1(R)$. Since $\AAA$ is closed under extensions in $C_b(\MMM(R))$, $\AAA$ is an exact category in the natural way. Notice that in $\AAA$, every admissible exact sequence is degree-wised split.\\ 
(2) For an integer $n$, let $\AAA^n$ be the full sub-category of
$n$-spherical complexes in $\AAA$. Remember that a complex $X$ is called {\it{$n$-spherical complex}} if $\Homo_k(X)=0$ whenever $k \ne n$.\\ Since $\AAA^n$ is closed under the extensions in $\AAA$, $\AAA^n$ can be considered as an exact category in the natural way.\\
(3) Let $q(\AAA)$ be subcategory of all quasi-isomorphisms in $\AAA$. Then $q(\AAA)$ is a class of saturated extensional weak equivalences in $\AAA$. 
\end{nt}

\begin{rem} [Truncation functors on $\AAA$] \label{ind fun 2} $ $
\\
(1) For any integer $c$, the functors 
$\tau_{\geqq c},\ \tau_{\leqq c-1}:C(\MMM(R)) \to C(\MMM(R))$
induces a functors 
$\tau_{\geqq c},\ \tau_{\leqq c-1}:\AAA \to \AAA$.

\begin{proof}
For any complex $X$ in $\AAA$, $\ker d_c^X$ and 
$\im d_c^X$ are finitely generated free $R$-modules, so $\tau_{\geqq c}X$
and $\tau_{\leqq c-1}X$ are bounded complexes of finitely generated
free $R$-modules. Considering the long exact sequences associated with the short exact sequence (\ref{short exact}) in \textbf{Notation \ref{gncc}}, we learn that their homology groups are torsion $R$-modules.  
\end{proof}

\noindent
(2) For any complex $X$ in $\AAA$ and an integer $n$, we have the following split exact sequence:
$$ 0 \to \tau_{\geqq n+1}X \to X \to \tau_{\leqq n}X \to 0.$$

\begin{proof}
For simplicity we will put $Y=\tau_{\geqq n+1}X$ and $Z=\tau_{\leqq n}X$ and let $f:Y \to X$ and $g:X \to Z$ be the natural morphisms. Since the sequence above is degree-wised split, there are  $s_{n}:Z_{n} \to X_{n}$ and $t_n:X_n \to Y_n$ such that
$t_n \circ f_n=\id$, $g_n \circ s_n=\id$, $t_n \circ s_n=0$, $g_n \circ f_n =0$, and $f_n \circ t_n + s_n \circ g_n =\id$.\\  
Obviously $g_k:X_k \to Z_k$ is an isomorphism for any $k<n$, and  $f_k:Y_k \to X_k$ is an isomorphism for any $k>n$.
We obtain  
$$  d_{n}^X \circ s_{n}=g_{n-1}^{-1} \circ g_{n-1} \circ d_n^X \circ s_n =
g_{n-1}^{-1} \circ d_n^Z \circ g_n \circ s_n=g_{n-1}^{-1} \circ d_n^Z\ \ .$$
$$ t_n \circ d_{n+1}^X = t_n \circ d_{n+1}^X \circ f_{n+1} \circ f_{n+1}^{-1}
= t_n \circ f_n \circ d_{n+1}^Y \circ f_{n+1}^{-1} =d_{n+1}^Y \circ f_{n+1}^{-1}\ \ .$$
This means that if we put 
$$ u_k=
\begin{cases}
0 & \text{if $k<n$}\\
t_n & \text{if $k=n$}\\
f_k^{-1} & \text{if $k>n$}
\end{cases}
,\ \ v_k= 
\begin{cases}
g_k^{-1} & \text{if $k<n$}\\
s_n & \text{if $k=n$}\\
0  & \text{if $k>n$}
\end{cases},$$
then $u :X \to Y$ and $v :Z \to X$ are morphisms of
chain complexes. By construction, we get $g \circ v =\id$, 
$u \circ f =\id$, $u \circ v =0$, $g \circ f=0$,
and $f \circ u + v \circ g =\id$.
Therefore the sequence $Y \overset{f}{\rightarrowtail} X \overset{g}
{\twoheadrightarrow} Z$ is a split exact sequence.\\
\end{proof}
\end{rem}

\begin{rem} $ $
\\
Since a complex in $\AAA$ is bounded, any morphism in $\AAA$ is
$n$-quasi-isomorphism for some integer $n$. Here for any integer $n$ a morphism of complexes $f:X \to Y$ in $\AAA$ is said to be {\it{$n$-quasi-isomorphism}},
if $\Homo_k(\Cone f)=0$ for any $k \leqq n$. This means $\Homo_k(f)$ is an 
isomorphism for $k<n$ and epimorphism for $k=n$.\\  
\end{rem}

\begin{prop} \label{Comp 1} $ $
\\
For any integer $n$, the inclusion functor $\AAA^n \to \AAA$ 
gives a homotopy equivalence 
$$S_{\bullet}(\AAA^n;q) \to S_{\bullet}(\AAA;q).$$
\end{prop}

\begin{proof}
We will use the cell filtrations theorem in \cite{Wal85} p.361 Theorem 1.7.1. To do so, we shall only verify that any morphism in $\AAA$ admits a {\it{cellular factorization}}. Here a cellular factorization of $n$-quasi isomorphism 
$f:X \to Y$ in $\AAA$ means a factorization of $f$:
$$X^n \overset{i^n}{\rightarrowtail} X^{n+1} \overset{i^{n+1}}{\rightarrowtail} \cdots \overset{i^{m-1}}{\rightarrowtail}
X^m \overset{g}{\to}Y$$
such that \\
(a) each $i^n$ is a degree-wised split monomorphism,\\
(b) $g$ is a quasi-isomorphism, and\\
(c) each $X^k/X^{k-1}$ is in $\AAA^k$\\
\\
\noindent
First we will prove the following \textbf{Claim 1}.\\
\\
\textbf{Claim 1}\\
Let $f:X \to Y$ be $n$-quasi-isomorphism in $\AAA$.
Then there are morphisms of complexes $g:X \to Z$ and $h:Z \to Y$ such that\\
(1) $f=h \circ g$,\\
(2) $\Cone h$ is homotopy equivalent to $\tau_{\geqq n+2}\Cone f$ , and\\
(3) $\Cone g$ is $(n+1)$-spherical complex.\\

\begin{proof}[Proof of {\rm \textbf{Claim 1}}]
By \textbf{Remark \ref{ind fun 2}} (2), we know $W=\tau_{\geqq n+2}\Cone f$ is a direct summand of $\Cone f$. Hence we have a canonical morphism
$a:\Cone f \to W.$
We put 
$$Z:=\Cone(Y \overset{\scriptstyle{j_2(f)}}{\to} \Cone f \overset{\scriptstyle{a}}{\to} W)[1]$$
(For the definition of $j_2(f)$, see \textbf{Notation \ref{gncc}} (5)) and
$h:Z \to Y$ is a canonical morphism. Then we have the commutative diagram:
$$\xymatrix{
X \ar[r]^f \ar@{-->}[d]_{g}& Y \ar[r]^{b} \ar[d]^{\id_Y} & \Cone f \ar[d]^{a} \\
Z \ar[r]^h & Y \ar[r] & \tau_{\geqq n+2}\Cone f & .
}$$

\noindent
{\it{The definition of $g$}}\\
We take a homotopy $H$ between $0$ and $abf$, that is a family of morphisms
$H_n:X_n \to W_{n+1}$ such that 
$H_{n-1} d_n^X + d_{n+1}^{\Cone f}H^n=-a_n b_n f_n.$ 
Put $g_n=\begin{pmatrix} H_n \\ f_n \end{pmatrix}$. We can easily check
that this map makes the diagram above commutative.\\
\\
\noindent
By the octahedral axiom, there is a distinguished triangle
$$\Cone g \to \Cone f \overset{a}{\to} \tau_{\geqq n+2}\Cone f \overset{+1}{\to}$$ 
in the homotopy category of $\AAA$. Using long exact sequence associated with the distinguished triangle above,
we learn that $\Cone g$ is a $(n+1)$-spherical complex.
\end{proof}

\noindent
Next we will prove the following \textbf{Claim 2}.\\
\\
\textbf{Claim 2}\\
In \textbf{Claim 1}, we can take $g$ as a degree-wised split monomorphism.

\begin{proof}[Proof of \textbf{Claim 2}]
Replace $Z$ with $\Cyl g$, $g$ with $j_1(g)$, and $h$ with $h \circ p(g)$. 
\end{proof}  

\noindent
Since any complex in $\AAA$ is bounded, Using \textbf{Claim 2} and
induction, we arrived at the result. 
\end{proof}

\begin{prop}
For any integer $k$, $\tau_{\leqq k}, \tau_{\geqq k}:\AAA^n \to
\AAA^n$ are exact functors. 
\end{prop}

\begin{proof}
By the following more abstract statement \textbf{Lemma \ref{snake lemma exer}}, this proposition is easily proved. (Notice also \textbf{Proposition \ref{exact check}}.)
\end{proof}

\begin{lem} \label{snake lemma exer} $ $
\\
Let $\XXX$ be an abelian category and 
$$0 \to X \to Y \to Z \to 0$$
a short exact sequence of complexes in $C(\XXX)$. Suppose
$\Homo_{n-1}(X)=0$ or $\Homo_n(Z)=0$ ,then the following sequences are
exact{\rm :}
\begin{equation} \label{short exact2}
0 \to \ker d_n^{X} \to \ker d_n^{Y} \to \ker d_n^{Z} \to 0
\end{equation}
\begin{equation} \label{short exact3}
0 \to \im d_n^{X} \to \im d_n^{Y} \to \im d_n^{Z} \to 0
\end{equation}
\end{lem}

\begin{proof}
By virtue of considering the diagram below, once we prove either 
(\ref{short exact2}) or (\ref{short exact3}) is exact, we obtain exactness of the other sequence by the $3 \times 3$ lemma.
$$\xymatrix{
& 0 \ar[d] & 0 \ar[d] &  0 \ar[d]  \\
0 \ar[r] & \ker d_n^{X} \ar[r] \ar[d] & \ker d_n^{Y} \ar[r] \ar[d] & \ker d_n^{Z} \ar[d] & & \\
0 \ar[r] & X_n \ar[r] \ar[d] & Y_n \ar[r] \ar[d] & Z_n \ar[r] \ar[d] & 0\\
& \im d_n^{X} \ar[r] \ar[d] & \im d_n^{Y} \ar[r] \ar[d] & \im d_n^{Z} \ar[d]  \\
& 0 & 0 & 0 &,
}$$
where vertical sequences and first and middle line are exact sequences and
third lines are complexes.\\
\\
\textbf{Claim 1}\\
If $\Homo_{n-1}(X)=0$, then (\ref{short exact2}) is exact.

\begin{proof}[Proof of \textbf{Claim 1}] 
Applying the snake lemma to the following diagram 
$$\xymatrix{
0 \ar[r] & X_n \ar[r] \ar[d] & Y_n \ar[r] \ar[d] & Z_n \ar[r] \ar[d] & 0 & \text{(exact)}\\
0 \ar[r] & \ker d_{n-1}^X \ar[r] & \ker d_{n-1}^Y \ar[r] & \ker d_{n-1}^Z & &\text{(exact),}
}$$
we get an exact sequence
$$0  \to \ker d_n^X \to \ker d_n^Y \to \ker d_n^Z \to \Homo_{n-1}(X)=0 \text{ (exact)}.$$
\end{proof}  

\noindent
\textbf{Claim2}\\
If $\Homo_n(Z)=0$, then (\ref{short exact3}) is exact.

\begin{proof}[Proof of \textbf{Claim 2}]
First we note that $\ker d_n^Y \to \ker d_n^Z$ is an epimorphism, because
we know $\ker d_n^Z =\im d_{n+1}^Z$ by assumption and composition
$\im d_{n+1}^Y \to \ker d_n^Y \to \ker d_n^Z=\im d_{n+1}^Z$ is an epimorphism. Applying the snake lemma to the following diagram
$$\xymatrix{
0 \ar[r] & \ker d_n^X \ar[r] \ar[d] & \ker d_n^Y \ar[r] \ar[d] & \ker d_n^Z \ar[r] \ar[d] & 0 & \text{(exact)}\\
0 \ar[r] & X_n \ar[r] & Y_n \ar[r] & Z_n \ar[r] & 0 & \text{(exact),}
}$$
we get the exact sequence of (\ref{short exact3}).
\end{proof}

\noindent
Combining two Claims above, we get the result.
\end{proof}

\begin{cor} \label{Comp 2} $ $
\\
The inclusion functor $\lambda:\Kos^1(R) \to \AAA^0$ induces a homotopy equivalence
$$S_{\bullet}(\Kos^1(R);q) \isoto S_{\bullet}(\AAA^0;q).$$
\end{cor}

\begin{proof}
We put an exact functor $\kappa=\tau_{\leqq 0}\tau_{\geqq 0}:\AAA^0
\to \Kos^1(R)$. Obviously we know $\kappa \lambda =\id$.\\
We have the canonical natural transformations $\id \overset{u}{\leftarrow}
\tau_{\geqq 0} \overset{v}{\to} \lambda \kappa$ 
and for each $X \in \AAA^0$,
$u(X)$ and $v(X)$ are quasi-isomorphisms. Hence by \cite{Wal85} p.330 Proposition 1.3.1., we learn that $\kappa$ induces an inverse map of $S_{\bullet}(\Kos^1(R);q) \to S_{\bullet}(\AAA^0;q)$.  
\end{proof}

\begin{rem} [Well-known fact]  \label{wel fac2}$ $
\\
(c.f. \cite{TT90} 5.7) For any integer $n$, we have an isomorphism
$$K_n(\AAA;q) \isoto K_n(\MMM^1(R))$$ 
such that composition $K_n(\Kos^1(R);q) \to K_n(\AAA;q) \to K_n(\MMM^1(R))$ coincides with $K_n(\Homo_0)$. 
\end{rem}

\begin{proof}
Let $\BBB$ be the full subcategory of perfect complexes on $\Spec R$ which is quasi isomorphic to an object in $\AAA$. Then we have the following (non commutative) diagram of Waldhausen categories:
$$\xymatrix{
(\Kos^1(R);q) \ar[r]^{\text{I}} \ar[d]_{\Homo_0} \rtwocell \omit {<1.5>\ x} & (\AAA;q) \ar[d]^{\text{II}}\\
\MMM^1(R) \ar[r]_{\text{III}} & (\BBB;q)
}$$ 
where the morphisms I, II and III are canonical inclusions and the natural transformaion $x$ is defined by using the canonical surjection as follows.\\
$$\begin{bmatrix}
\xymatrix{
X_1 \ar[d]_{\scriptstyle{d^X}}\\
X_0
}
\end{bmatrix} 
\begin{matrix}
\xymatrix{
\to \ar@{}[d]\\
\to\\
}
\end{matrix} 
\begin{bmatrix}
\xymatrix{
0 \ar[d] \\
\Homo_0(X)
}
\end{bmatrix}$$
for any complex $X$ in $\Kos^1(R)$. Now we learn the diagram above induces commutative diagram on their higher algebraic $K$-groups by \cite{Wal85} p.330 Proposition 1.3.1 and the morphism I and II induces isomorphisms their higher algebraic $K$-groups. For I, see for example, the proof of the claim in the proof of Proposition 1.2. in \cite{Moc07}. For II, it follows by the approximation theorem \cite{Wal85}, \cite{TT90}, \cite{Sch06}. Hence we get the desired result.  
\end{proof}

\begin{proof}[Proof of Theorem \ref{comparison}]
It follows easily by \textbf{Proposition \ref{Comp 1}}, \textbf{Corollary \ref{Comp 2}} and \textbf{Remark \ref{wel fac2}}.
\end{proof}

\section{Semi-direct product of exact categories}

In this section, we will prove \textbf{Theorem \ref{split fib}}.

\begin{nt} [Semi-direct product of exact categories] \label{sdp} $ $
\\
Let $\FFF$ be $\MMM(R)$ or $\PPP(R)$. We will define the full sub-category
$\MMM^1(R) \ltimes \FFF$ of those complexes 
$X=[X_1 \overset{d^X}{\to} X_0]$ in $C_b(\MMM(R))$
such that  $X_i$ are in $\FFF$ for $i=0, 1$ and 
its boundary morphism $d^X$ is injective and its homology group are in 
$\MMM^1(R)$.\\ 
Clearly $\Kos^1(R)=\MMM^1(R) \ltimes \PPP(R)$ and we put $\CCC:=\MMM^1(R) \ltimes \MMM(R)$.\\
We say a sequence $X \to Y \to Z$ in $\MMM^1(R) \ltimes \FFF$
is an admissible short exact sequence if it is a short exact sequence
in $C_b(\MMM(R))$. Since $\MMM^1(R) \ltimes \FFF$ is closed under extensions in $C_b(\MMM(R))$, we learn that it can be considered as an exact category in the natural way. Considering the long exact sequence, we also learn 
$$\Homo_0:\MMM^1(R) \ltimes \FFF \to \MMM^1(R)$$ 
is an exact functor.
\end{nt}

\begin{lem} \label{comparison2} $ $
\\
The canonical inclusion functor 
$\Kos^1(R) \hookrightarrow \CCC$ induces a homotopy equivalence 
$$Q\Kos^1(R) \isoto Q\CCC.$$
\end{lem}

\begin{proof}
We will use Quillen's resolution theorem. The fact that $\Kos^1(R)$ is closed under extensions in $\CCC$ has already been remarked. Since a submodule of a finitely generated free $R$-module is also free, the fact that $\Kos^1(R)$ is also closed under taking kernel is easily proved. The fact that for any $Z$ in $\CCC$, there is an admissible epimorphism from object $Y$ in $\Kos^1(R)$ is proved as follows. First we can find surjective $Y_0 \to Z_0$ such that $Y_0$ is finitely generated free $R$-module.
Next put $Y_1=\ker(Y_0 \to Z_0 \to \Homo_0(Z))$, then diagram chases argument
shows there is a morphism $Y \to Z$ and this is what we wanted. 
\end{proof}

\begin{lem} \label{pre split fib} $ $
\\
For any non-negative integer $n$, we have a split exact sequence
$$0 \to K_n(\CCC^q) \to K_n(\CCC) \overset{K_n(\Homo_0)}{\to} K_n(\MMM^1(R)) \to 0$$ 
where $\CCC^q$ is the full subcategory of acyclic complexes in $\CCC$.
\end{lem}

\begin{proof}
As $\CCC^q$ is closed under extensions in $\CCC$, $\CCC^q$ can be considered as an exact category in the natural way. There is a category equivalence as exact categories 
$$\CCC \isoto 
\E(\CCC^q, \CCC,\MMM^1(R))\ \ \ \ ,$$
$$X \mapsto 
\begin{pmatrix}
\begin{bmatrix}
\xymatrix{
X_1 \ar[d]_{\scriptstyle{\id_{X_1}}}\\
X_1
}
\end{bmatrix} &
\begin{matrix}
\xymatrix{
\overset{\scriptstyle{\id_{X_1}}}{\to} \ar@{}[d]\\
\underset{\scriptstyle{d^X}}{\to}\\
}
\end{matrix} &
\begin{bmatrix}
\xymatrix{
X_1 \ar[d]_{\scriptstyle{d^X}}\\
X_0
}
\end{bmatrix} &
\begin{matrix}
\xymatrix{
\to \ar@{}[d]\\
\to\\
}
\end{matrix} &
\begin{bmatrix}
\xymatrix{
0 \ar[d] \\
\Homo_0(X)
}
\end{bmatrix}
\end{pmatrix}\ \ \ .$$
where $\E(\CCC^q,\CCC,\MMM^1(R))$ is the exact category of admissible exact sequences $X \to Y \to Z$ in $\CCC$ such that $X$ is in $\CCC^q$ and $Z$ is in $\MMM^1(R)$. Hence by the additivity theorem, we get the result.
\end{proof}

\begin{proof}[Proof of Theorem \ref{split fib}]
For any non-negative integer $n$, we have the following  commutative diagram: 
$$\xymatrix{
0 \ar[r] & K_n(\Kos^1(R)^q) \ar[r] \ar[d]_{\text{I}} & K_n(\Kos^1(R)) \ar[d]_{\text{II}} \ar[r] & K_n(\Kos^1(R);\qis) \ar[r] \ar[d]^{K_n(\Homo_0)} & 0\\
0 \ar[r] & K_n(\CCC^q) \ar[r] & K_n(\CCC) \ar[r]_{K_n(\Homo_0)} & K_n(\MMM^1(R)) \ar[r] & 0.
}$$
Here $K_n(\Homo_0)$ and the morphism II are isomorphisms by \textbf{Theorem \ref{comparison}} and \textbf{Lemma \ref{comparison2}}, respectively.
Now the morphism I is also an isomorphism. For
there are category equivalences
$$\Kos^1(R)^q \ni X \mapsto X_1 \in \PPP(R),$$ 
$$\CCC^q \ni X \mapsto X_1 \in \MMM(R).$$ 
and by the resolution theorem again. Since the bottom line is a split exact sequence by \textbf{Lemma \ref{pre split fib}}, the top line is also.
\end{proof}

\appendix

\section{Exercises in elementary divisors theory}

In this section, we will give an another proof of \textbf{Theorem \ref{split fib}}. We will start the following examples. From now on we will assume $R$ is commutative.

\begin{ex}[Schlichting's works] $ $
\\
In the notation \textbf{Theorem \ref{split fib}}, if we know the existence of the following fibration sequence
\begin{equation} \label{Schlichting}
K(\Kos^1(R)^q;\isom) \to K(\Kos^1(R);\isom) \to K(\Kos^1(R);\qis), 
\end{equation}
then we can easily prove \textbf{Theorem \ref{split fib}}. For the inclusion $\Kos^1(R)^q \hookrightarrow \Kos^1(R)$ has an obvious retraction up to the natural equivalence defined by
$$\Kos^1(R) \ni 
\begin{bmatrix}
\xymatrix{
X_0 \ar[d]_{\scriptstyle{d^X}}\\
X_1
}
\end{bmatrix}
\mapsto
\begin{bmatrix}
\xymatrix{
X_1 \ar[d]_{\scriptstyle{\id_{X_1}}}\\
X_1
}
\end{bmatrix} 
\in \Kos^1(R)^q.$$
(1) The first idea to get the fibration sequence above is using \lq\lq fibration theorem for Waldhausen categories with factorization\rq\rq \ (see \cite{Sch06} p.127 Theorem 11). But $\Kos^1(R)$ is not a Waldhausen category with factorization.

\begin{proof}
If a morphism $X \to 0$ from a non acyclic complex $X$ to the zero complex in $\Kos^1(R)$ has a factorization $X \overset{i}{\rightarrowtail} Z \overset{p}{\to} 0$ where $i$ is an admissible monomorphism and $p$ is a quasi isomorphism, then taking $\Homo_0$, we get the admissible monomorphism $\Homo_0(X) \rightarrowtail 0$. It is a contradiction. 
\end{proof} 
(2) The second idea is using \lq\lq homotopy fibration of $K$-theory spaces associated to Karoubi-Schlichting filtration\rq\rq \ (see \cite{Sch04} p.1097 Theorem 2.1). We will prove in \textbf{Proposition \ref{exercise}} that $\Kos^1(R)^q$ is idempotent complete right s-filtering in $\Kos^1(R)$ in the sense of \cite{Sch04} p.1091 Definition 1.3, 1.5. Then we can easily prove that $\Homo_0$ induces a category equivalence 
\begin{equation} \label{quotient equivalent}
\Kos^1(R)/\Kos^1(R)^q \isoto \MMM^1(R)
\end{equation}
by the universal property of the quotient exact category. Now there is the following commutative diagram of spaces:
$$\xymatrix{
K(\Kos^1(R)^q) \ar[r] \ar@{=}[dd] & K(\Kos^1(R)) \ar[r] \ar@{=}[dd] & K(\Kos^1(R);q) \ar[d]^{K(\Homo_0)}\\
& & K(\MMM^1(R))\\
K(\Kos^1(R)^q) \ar[r] & K(\Kos^1(R)) \ar[r] & K(\Kos^1(R)/\Kos^1(R)^q). \ar[u]\\
}$$
Here the vertical morphisms are isomorphisms by \textbf{Theorem \ref{comparison}} and the equivalence (\ref{quotient equivalent}), and the bottom line is a fibration sequence by Schlichting fibration theorem above. Hence we learn the existence of (\ref{Schlichting}). So we have an another proof of \textbf{Theorem \ref{split fib}}. 
\end{ex}

\noindent
So main purpose of this section is to prove the following proposition.

\begin{prop} \label{exercise} $ $
\\
$\Kos(R)^q$ is idempotent complete right s-filtering in $\Kos^1(R)$.
\end{prop}

\begin{proof}
First we list up what we shall prove as follows. (see \cite{Sch04} Ibid.)\\
(a) $\Kos(R)^q$ is idempotent complete.\\
(b) For an admissible exact sequence
\begin{equation} \label{assertion b}
X \rightarrowtail Y \twoheadrightarrow Z
\end{equation}
in $\Kos^1(R)$, $X$ and $Z$ are acyclic if and only if $Y$ is.\\
(c) Every morphism $X \to Y$ from an acyclic complex in $\Kos^1(R)$ to a complex in $\Kos^1(R)$ factor through an acyclic complex $Z$ in $\Kos^1(R)$ such that the morphism $X \twoheadrightarrow Z$ is an admissible epimorphism.\\
(d) For every admissible monomorphism $X \rightarrowtail Y$ from an acyclic complex in $\Kos^1(R)$ to a complex in $\Kos^1(R)$, there is an admissible epimorphism $Y \twoheadrightarrow Z$ such that the composition $X \rightarrowtail Z$ is an admissible monomorphism in $\Kos^1(R)$. (By the assertion (b), it is also the admissible monomorphism in $\Kos^1(R)^q$.)\\
\\
\noindent
The assertion (b) is easily checked by considering the long exact sequence associated with (\ref{assertion b}). The assertion (a) and (c) are follow from the following \textbf{Claim}.\\ 
\\
\noindent
\textbf{Claim}\\
For every morphism $f:X \to Y$ from an acyclic complex in $\Kos^1(R)^q$ to an complex in $\Kos^1(R)$, $\im f$ is an acyclic complex and the canonical morphism $X \to \im f$ is an admissible epimorphism in $\Kos^1(R)^q$.\\

\begin{proof}[Proof of \textbf{Claim}]
First notice that for any morphism $f:X \to Y$ in $\Kos^1(R)$, we have the following sequence \begin{equation} \label{claim}
0 \to \ker f \to X \to \im f \to 0
\end{equation}
in $C_b(\MMM(R))$. Since a submodule of a finitely generated free $R$-module is also free, we learn that every components in $\ker f$ and $\im f$ are free. Considering the long exact sequence associated with (\ref{claim}), we learn $\ker f$ and $\im f$ are objects in $\Kos^1(R)$. If $X$ is acyclic, then by applying the assertion (b) for (\ref{claim}), we learn that $\ker f$ and $\im f$ are also acyclic, and therefore the short exact sequence (\ref{claim}) is admissible in $\Kos^1(R)^q$.
\end{proof} 
\noindent
From now on, we will prove the assertion (d). We will fix the notations. Let $i:X \rightarrowtail Y$ be an admissible monomorphism from an acyclic complex in $\Kos^1(R)$ to a complex in $\Kos^1(R)$. So there is an admissible exact sequence
in $\Kos^1(R)$.
$$X \overset{i}{\rightarrowtail} Y \overset{p}{\twoheadrightarrow} W.$$
By elementary divisors theory, we learn that $Y$ and $W$ can be identified with the following type complexes.
$Y=V \oplus X \oplus U$ and $W=V \oplus U$. Here $U=[U \overset{\id}{\to} U]$, $V=[R^{\oplus n} \overset{d^V}{\to} R^{\oplus n}]$ and
$$d^V=
\begin{pmatrix}
a_1 & & & &\\
& a_2 & & &\\
& & a_3 & &\\
& & & \cdots &\\
& & & & a_n\\
\end{pmatrix}$$
where $a_1||a_2||\cdots ||a_r$ are non units and $n$ is rank of $\Homo_0(Y)$. So we can write $i$ and $p$ as the matrix descriptions:
$$i_s=\begin{pmatrix}i_1^s\\i_2^s\\i_3^s\end{pmatrix}:
X_s \to V_s \oplus X_s \oplus U,$$
$$p_s=\begin{pmatrix}p_{11}^s & p_{12}^s & p_{13}^s\\ p_{21}^s & p_{22}^s & p_{23}^s \end{pmatrix}:V_s \oplus X_s \oplus U \to V_s \oplus U$$  
for $s=0,1$. Since $i$ is an admissible monomorphism and every admissible exact sequence in $\Kos^1(R)$ is degree-wised split, we have an morphism
$$h=\begin{pmatrix} h_1 & h_2 & h_3 \end{pmatrix}:V_0 \oplus X_0 \oplus U \to X_0$$such that $hi_0=\id$. Then we put $Z=X \oplus U$ and define $q:Y \to Z$ as follows.
$$q_0=\begin{pmatrix} h_1 & h_2 & h_3\\ p_{21}^0 & p_{22}^0 & p_{32}^0 \end{pmatrix}:V_0 \oplus X_0 \oplus U \to X_0 \oplus U,$$
$$q_1=\begin{pmatrix} {(d^X)}^{-1}h_1 d^V & {(d^X)}^{-1}h_2d^X & {(d^X)}^{-1}h_3\\ p_{21}^0d^V & p_{22}^0d^X & p_{32}^0 \end{pmatrix}:V_1 \oplus X_1 \oplus U \to X_1 \oplus U.$$
Now we can easily check that $q$ is a morphism of chain complex and $qi=\begin{pmatrix} \id\\0\end{pmatrix}$. Next we will prove $q$ is an admissible epimorphism in $\Kos^1(R)$. Since $\Kos^1(R)$ is closed under taking kernel (this is proved in the poof of \textbf{Claim}), we shall only prove that $q$ has degree-wised sections. To do so first we pick up a morphism 
$$l:=\begin{pmatrix}l_{11} & l_{12}\\ l_{21} & l_{22}\\ l_{31} & l_{32} \end{pmatrix}:V_0 \oplus U \to V_0 \oplus X_0 \oplus U$$
such that it is a section of $p_0$. Then we can easily check the identity
$$\begin{pmatrix}
h_1 & h_2 & h_3\\
p_{21}^0 & p_{22}^0 & p_{23}^0
\end{pmatrix}
\begin{pmatrix}
i_1^0 & l_{12}\\
i_2^0 & l_{22}\\
i_3^0 & l_{32}
\end{pmatrix}  
{\begin{pmatrix}
\id & h_1l_{12}+h_2l_{22}+h_3l_{32}\\
0 & \id
\end{pmatrix}}^{-1}
=\id$$
Hence we learn that $q_0$ possesses a section. Similarly we learn that $q_1$ also does. Hence we complete proving the assertion (d).
\end{proof}

\begin{center}{\it
Satoshi Mochizuki\\
Graduate school of Mathematical Sciences,
The University of Tokyo, 3-8-1 Komaba, Meguro-ku
Tokyo 153-8914, JAPAN}
\\
\tt{E-mail:mochi@ms.u-tokyo.ac.jp}
\end{center}


\begin{thebibliography}{99}
\bibitem[Gab62]{Gab62}
P. Gabriel, {\it{Des cat\'egories ab\'eli\`ennes}}, Bull. Soc. Math. France \textbf{90}, (1962), p.323-448.
\bibitem[Kel90]{Kel90}
B. Keller, {\it{Chain complexes and stable categories}}, manus. math. \textbf{67} (1990), p.379-417.
\bibitem[Moc07]{Moc07}
S. Mochizuki, {\it{Gersten conjecture for commutative discrete valuation rings}}, avivable at http://www.math.uiuc.edu/K-theory/0819 (2007).
\bibitem[Qui73]{Qui73} 
D. Quillen, {\it{Higher algebraic K-theory I}}, In Higher K-theories, Springer Lect. Notes Math. \textbf{341} (1973), p.85-147.
\bibitem[Sch04]{Sch04}
M. Schlichting, {\it{Delooping the $K$-theory of exact categories}}, Topology \textbf{43} (2004), p.1089-1103.
\bibitem[Sch06]{Sch06}
M. Schlichting, {\it{Negative $K$-theory of derived categories}}, Math. Z. \textbf{253} (2006), p.97-134.
\bibitem[TT90]{TT90}
R. W. Thomason, T. Trobaugh, {\it{Higher K-theory of schemes and of derived categories}}, In The Grothendieck Festscrift,Vol III, (1990), p.247-435. 
\bibitem[Wal78]{Wal78}
F. Waldhausen, {\it{Algebraic K-theory of generalized free products}}, Ann. of Math. \textbf{108} (1978), p.135-256.
\bibitem[Wal85]{Wal85}
F. Waldhausen, {\it{Algebraic K-theory of spaces}}, In Algebraic and geometric topology, Springer Lect. Notes Math \textbf{1126} (1985), p.318-419.
\end{thebibliography}
\end{document}